\documentclass[12pt]{article}

\usepackage{amssymb}
\usepackage{amsmath,amsthm}
\usepackage[dvips]{graphicx}
\usepackage{wrapfig}
\usepackage{bm}

\newtheorem{theorem}{Theorem}[section]

\newtheorem{lemma}[theorem]{Lemma}

\theoremstyle{definition}

\theoremstyle{remark}

\title{A note on unavoidable sets for a spherical curve of reductivity four}

\author{Kenji Kashiwabara , Ayaka Shimizu }
\date{\today}

\begin{document}

\maketitle

\begin{abstract}
The reductivity of a spherical curve is the minimal number of a local transformation called an inverse-half-twisted splice required to obtain a reducible spherical curve from the spherical curve. 
It is unknown if there exists a spherical curve whose reductivity is four. 
In this paper, an unavoidable set of configurations for a spherical curve with reductivity four is given by focusing on 5-gons. 
It has also been unknown if there exists a reduced spherical curve which has no 2-gons and 3-gons of type A, B and C. 
This paper gives the answer to the question by constructing such a spherical curve. 
\end{abstract}



\section{Introduction}

\begin{figure}[ht]
\begin{center}
\includegraphics[width=80mm]{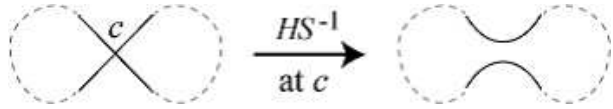}
\caption{An inverse-half-twisted splice operation at a crossing $c$. Broken curves represent the outer connections. }
\label{hs}
\end{center}
\end{figure}

A {\it spherical curve} is a closed curve on $S^2$, where self-intersections, called {\it crossings}, are double points intersecting transversely. 
In this paper, spherical curves are considered up to ambient isotopy of $S^2$, and two spherical curves which are transformed into each other by a reflection are assumed to be the same spherical curve. 
A spherical curve is {\it trivial} if it has no crossings. 
A spherical curve $P$ is {\it reducible} if one can draw a circle on $S^2$ which intersects $P$ transversely at just one crossing of $P$. 
Otherwise, it is said to be {\it reduced}. 
An {\it inverse-half-twisted splice}, denoted by ${HS}^{-1}$, at a crossing of a spherical curve $P$ is a splice on $P$ which yields another spherical curve (not a link projection) as shown in Figure \ref{hs}. 
An inverse-half-twisted splice does not preserve an orientation of a spherical curve. 
Hence ${HS}^{-1}$ is a different local transformation from the splice called a ``smoothing'' in knot theory. 
\cite{IS} shows that for every pair of two nontrivial reduced spherical curves $P$ and $P'$, there exists a finite sequence of ${HS}^{-1}$s and its inverses which transform $P$ into $P'$ such that a spherical curve at each step of the sequence is also reduced. 
This implies that all nontrivial reduced spherical curves are connected by ${HS}^{-1}$s and its inverses. 
The {\it reductivity} of a nontrivial spherical curve $P$ is defined to be the minimal number of inverse-half-twisted splices, ${HS}^{-1}$s, which are required to obtain a reducible spherical curve from $P$. 
The reductivity tells us how reduced a spherical curve is like the connectivity in graph theory. 
In \cite{S}, it is shown that every nontrivial spherical curve has the reductivity four or less. 
Also, in \cite{S} and \cite{OS}, it is mentioned that there are infinitely many spherical curves with reductivity 0, 1, 2 and 3. 
At the moment, the following problem is open: 

\phantom{x}
\noindent {\bf Problem A. (\cite{S})}. \ For any nontrivial spherical curve, is the reductivity three or less?\\
\phantom{x}

\begin{figure}[ht]
\begin{center}
\includegraphics[width=60mm]{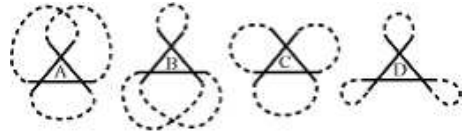}
\caption{The 3-gons of type A, B, C and D. Broken curves represent outer connections. }
\label{3-gons}
\end{center}
\end{figure}

\begin{figure}[ht]
\begin{center}
\includegraphics[width=140mm]{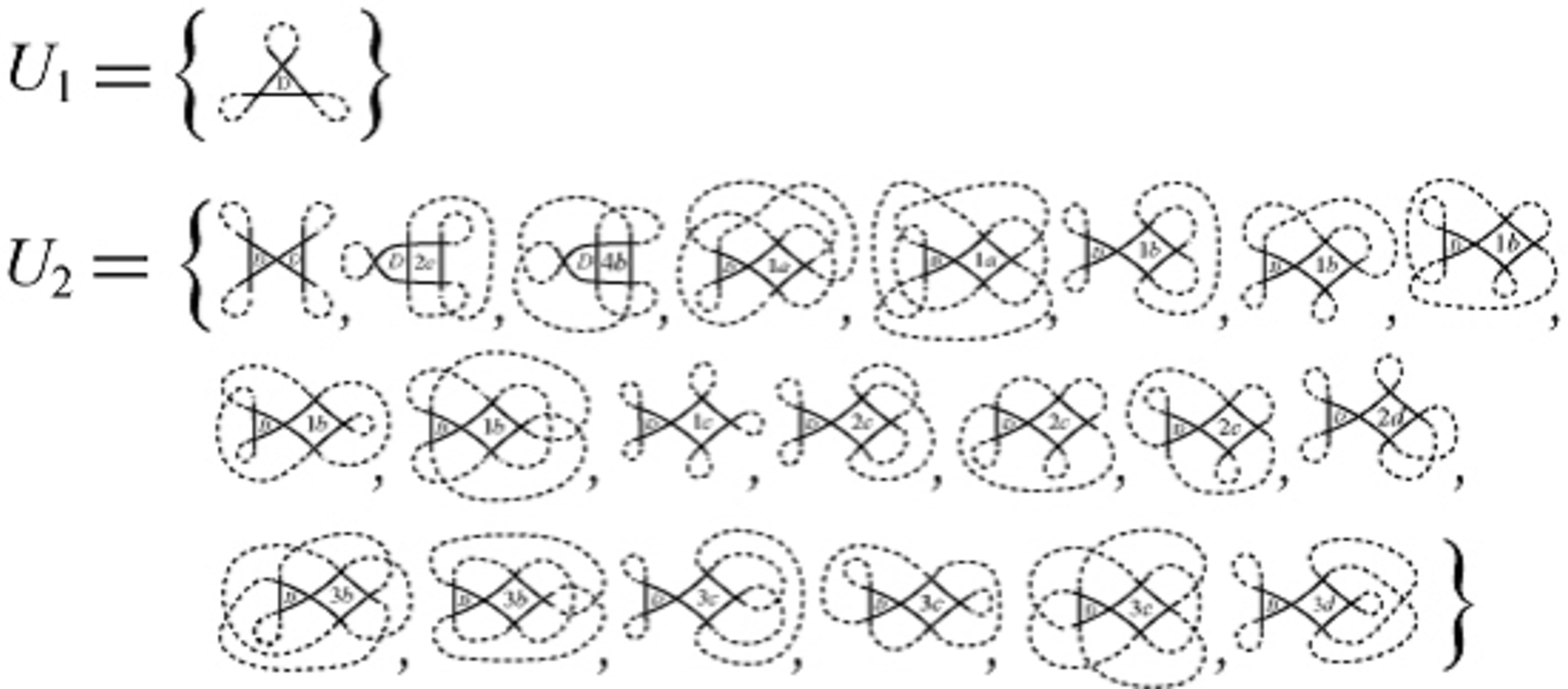}
\caption{Unavoidable sets with outer connections $U_1$ and $U_2$ for a spherical curve of reductivity four. Broken curves represent outer connections. }
\label{u-set-4}
\end{center}
\end{figure}

\noindent In other words, it is unknown if there exists a spherical curve whose reductivity is four. 
An {\it unavoidable set} of configurations for a spherical curve in a class is a set of configurations with the property that any spherical curve in the class has at least one member of the set (see, for example, \cite{CZ}). 
It is important to find unavoidable sets for a spherical curve of reductivity four from various viewpoints. 
In \cite{S}, 3-gons were classified into four types considering outer connections as shown in Figure \ref{3-gons} and the unavoidable set $U_1$, shown in Figure \ref{u-set-4}, of configurations with outer connections for a spherical curve with reductivity four was given. 
\begin{figure}[ht]
\begin{center}
\includegraphics[width=90mm]{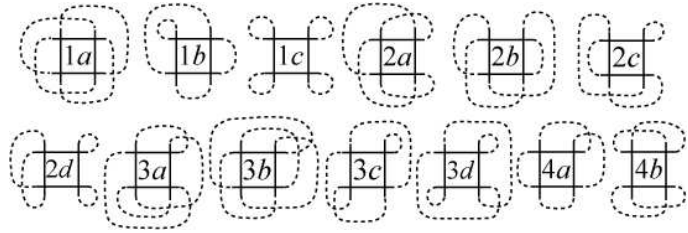}
\caption{The 13 types of 4-gons. }
\label{4-gons}
\end{center}
\end{figure}
\begin{figure}[ht]
\begin{center}
\includegraphics[width=80mm]{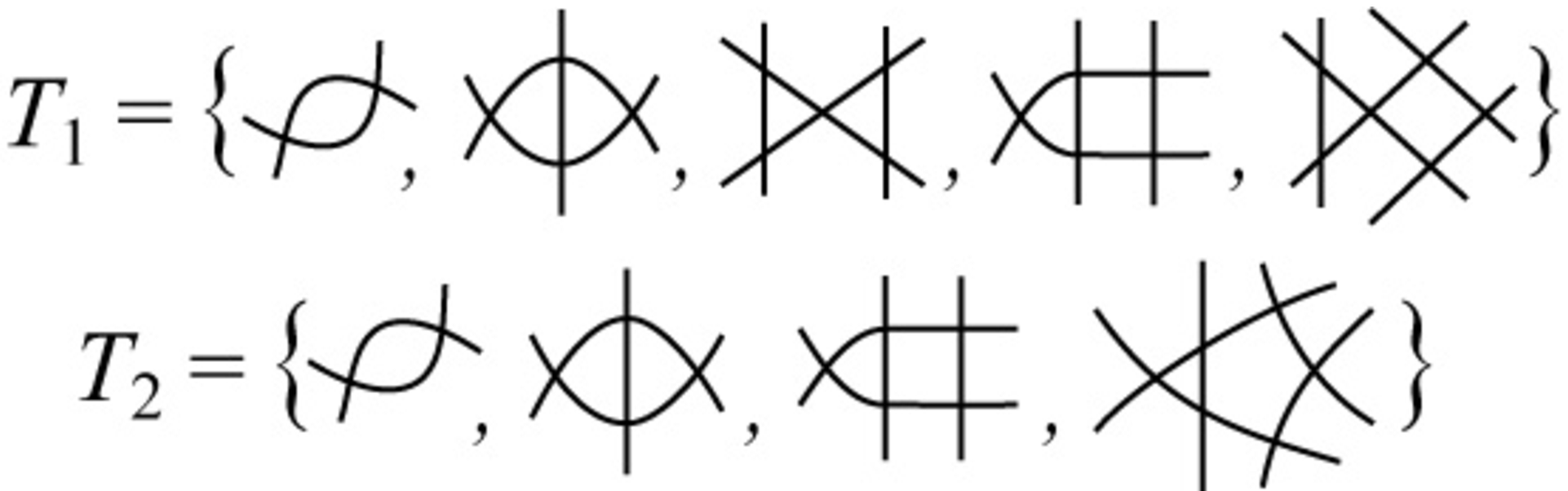}
\caption{Unavoidable sets of configurations (with any outer connections) $T_1$ and $T_2$ for a reduced spherical curve. }
\label{u-set-r}
\end{center}
\end{figure}
The unavoidable set $U_1$ was obtained by the following two facts; 
the first one is that every nontrivial reduced spherical curve has a 2-gon or 3-gon (\cite{AST}). 
The second one is that if a spherical curve has a 2-gon or a 3-gon of type A, B or C, then the reductivity is three or less (\cite{S}). 
The following problem was also posed in \cite{S}: 

\phantom{x}
\noindent {\bf Problem B. (\cite{S})} \ 
Is the set consisting of a 2-gon, 3-gons of type A, B and C an unavoidable set for a reduced spherical curve? \\
\phantom{x}

\noindent If the answer to Problem B is ``yes'', then the answer to Problem A is also ``yes''. 
However, the following theorem gives the negative answer to Problem B: 

\phantom{x}
\begin{theorem}
There exists a reduced spherical curve which has no 2-gons and 3-gons of type A, B and C. 
\label{d-thm}
\end{theorem}
\phantom{x}

\noindent (See Figures \ref{d-kps} and \ref{d-3} in Section 2.) 
In \cite{OS}, 4-gons were classified into 13 types as shown in Figure \ref{4-gons} and the unavoidable set $U_2$, in Figure \ref{u-set-4}, for a spherical curve with reductivity four was given by combining 3-gons and 4-gons based on an unavoidable set $T_1$ in Figure \ref{u-set-r} for a nontrivial reduced spherical curve which was obtained in \cite{S} in the same way to the four-color-theorem. 
Note that a necessary condition for a spherical curve with reductivity four was also given using the notion of the warping degree in \cite{KS}. 
In this paper, 5-gons are classified in a systematic way which can be used for general $n$-gons (in Section 3) and another unavoidable set for a spherical curve with reductivity four is given: 

\phantom{x}

\begin{theorem}
The set $U_3$ shown in Figure \ref{u-set} is an unavoidable set for a spherical curve with reductivity four. 
\begin{figure}[ht]
\begin{center}
\includegraphics[width=130mm]{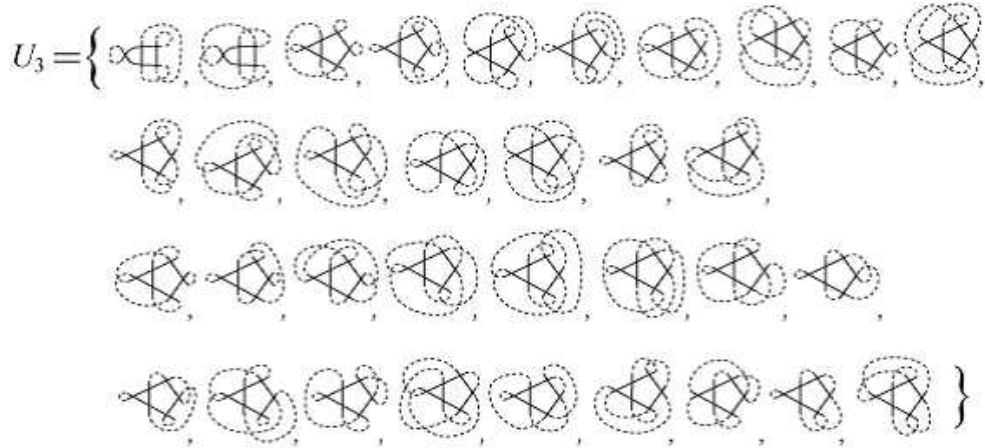}
\caption{An unavoidable set of configurations with outer connections for a spherical curve of reductivity four. }
\label{u-set}
\end{center}
\end{figure}
\label{thm-u}
\end{theorem}
\phantom{x}

\noindent Theorem \ref{thm-u} would be useful for constructing a spherical curve with reductivity four (or showing that there are no such spherical curves), or detecting the reductivity for spherical curves which have no 2-gons and 3-gons of type A, B and C. 
The rest of the paper is organized as follows: 
In Section 2, Theorem \ref{d-thm} is shown. 
In Section 3, 5-gons are classified into 56 types. 
In Section 4, Theorem \ref{thm-u} is proved. 
In Appendix, the 5-gons on chord diagrams are listed.

\section{Proof of Theorem \ref{d-thm}}

In this section, Theorem \ref{d-thm} is shown. \\

\begin{figure}[ht]
\begin{center}
\includegraphics[width=130mm]{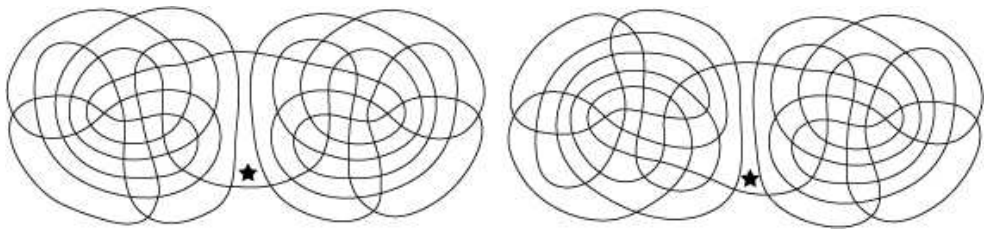}
\caption{Reduced spherical curves without 2-gons, 3-gons of type A, B and C. }
\label{d-kps}
\end{center}
\end{figure}

\phantom{x}
\noindent {\rm Proof of Theorem \ref{d-thm}.} \ 
\noindent The spherical curves depicted in Figure \ref{d-kps} are reduced, and have no 2-gons and 3-gons of type A, B and C. 
The point is that there are no 2-gons, and all the 3-gons are of type D. 
\hfill$\square$

\phantom{x}

\begin{figure}
\begin{center}
\includegraphics[width=130mm]{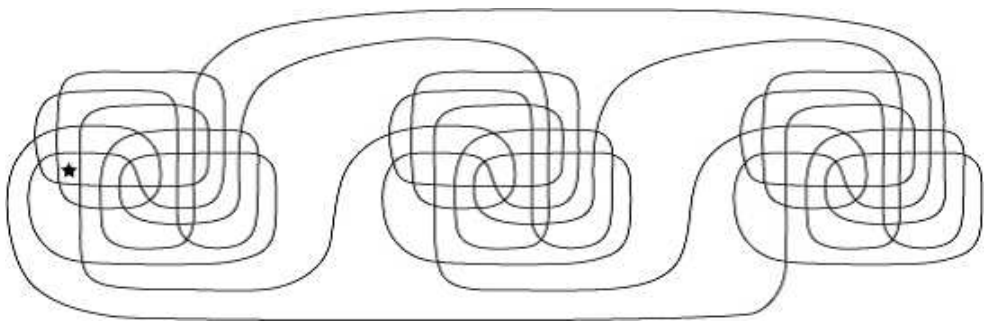}
\caption{Reduced spherical curve without 2-gons, 3-gons of type A, B and C. }
\label{d-3}
\end{center}
\end{figure}

\noindent Note that the spherical curves shown in Figure \ref{d-kps} have the reductivity one, not four,  because an inverse-half-twisted splice at a crossing at the middle 4-gons with a star derives a reducible spherical curve. 
Another example is shown in Figure \ref{d-3}. 
The reductivity of the spherical curve in Figure \ref{d-3} is not four because it has a 4-gon, with a star in the figure, of type 4a; 
it is shown in \cite{OS} that if a spherical curve has a 4-gon of type 4a, then the reductivity is three or less. 

\begin{figure}
\begin{center}
\includegraphics[width=70mm]{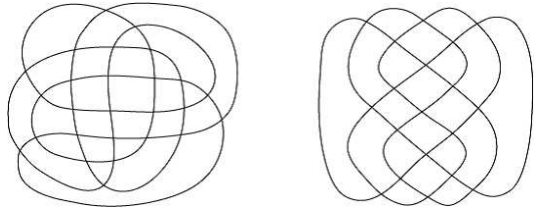}
\caption{Reduced spherical curves without 2-gons, 3-gons of type A and B. }
\label{cd-kps}
\end{center}
\end{figure}

In \cite{S}, a reduced spherical curve which has no 2-gons and 3-gons of type A and B was given. 
Further spherical curves are shown in Figure \ref{cd-kps}. 

\section{5-gons}

\begin{figure}[ht]
\begin{center}
\includegraphics[width=140mm]{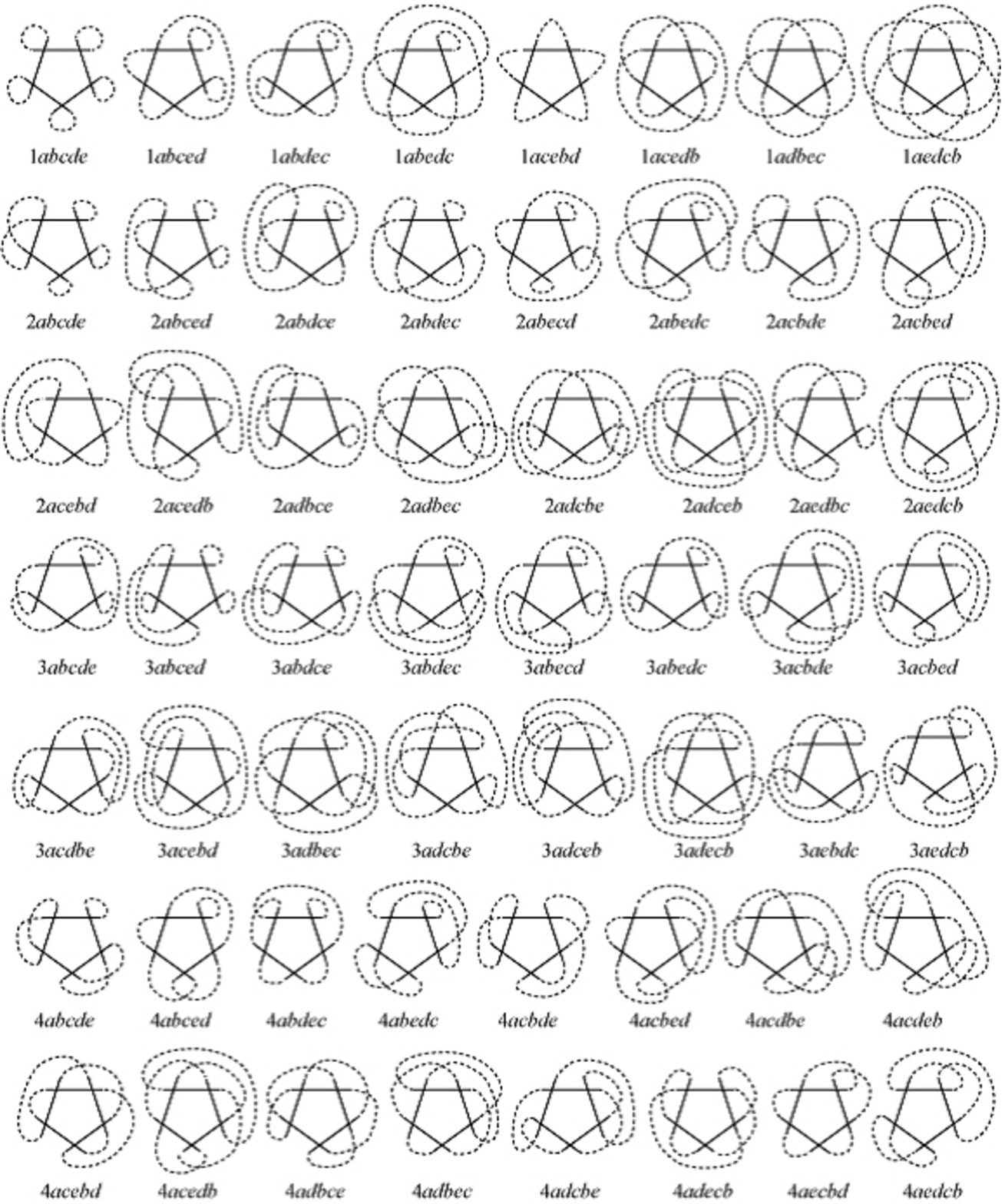}
\caption{All the 5-gons of a spherical curve with outer connections. }
\label{a-5-gons}
\end{center}
\end{figure}

In this section, the following lemma is shown: 

\phantom{x}
\begin{lemma}
5-gons of a spherical curve are divided into the 56 types in Figure \ref{a-5-gons} with respect to outer connections. 
\end{lemma}
\phantom{x}

\begin{figure}[ht]
\begin{center}
\includegraphics[width=70mm]{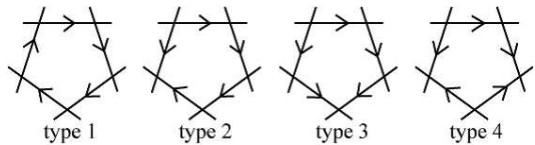}
\caption{5-gons of type 1 to 4 with relative orientations of the sides. }
\label{type1-4}
\end{center}
\end{figure}

\begin{figure}[ht]
\begin{center}
\includegraphics[width=120mm]{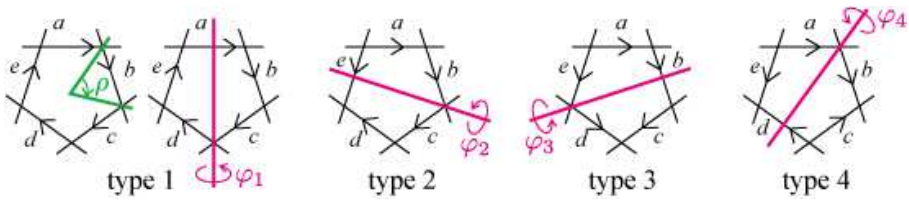}
\caption{The rotation and reflections on the 5-gons. (One of the relative orientations is shown by arrows at each type.) }
\label{type1-4-sym}
\end{center}
\end{figure}

\noindent {\rm Proof. } 
There are four types of 5-gons when relative orientations of the five sides are considered. 
The 5-gons of type 1 to 4 are illustrated in Figure \ref{type1-4}, where one of the relative orientations are shown by arrows. 
Let $a, b, c, d$ and $e$ be the sides of a 5-gon as illustrated in Figure \ref{type1-4-sym}. 
The 5-gon of type 1 has two types of symmetries: the ($2 \pi / 5$)-rotation $\rho$ and the reflection $\varphi _1$ defined by the following permutations
\begin{equation*}
\rho =
\begin{pmatrix}
a & b & c & d & e \\
b & c & d & e & a
\end{pmatrix}
, \ \varphi _1 =
\begin{pmatrix}
a & b & c & d & e \\
a & e & d & c & b
\end{pmatrix}
. 
\end{equation*}
The 5-gons of type 2, 3 and 4 have the reflection symmetries $\varphi _2$, $\varphi _3$ and $\varphi _4$ defined by the following permutations, respectively: 
\begin{equation*}
\varphi _2 =
\begin{pmatrix}
a & b & c & d & e \\
d & c & b & a & e
\end{pmatrix}
, \ \varphi _3 =
\begin{pmatrix}
a & b & c & d & e \\
c & b & a & e & d
\end{pmatrix}
, \ \varphi _4 =
\begin{pmatrix}
a & b & c & d & e \\
b & a & e & d & c
\end{pmatrix}
. 
\end{equation*}

Now let a 5-gon be a part of a spherical curve on $S^2$. 
Let $a, b, c, d$ and $e$ be sides of the 5-gon located as same as Figure \ref{type1-4-sym}. 
Fix the orientation of $a$ as $e$ to $b$. 
By reading the sides up as one passes the spherical curve, a cyclic sequence consisting of $a, b, c, d$ and $e$ is obtained. 
In particular, a sequence starting with $a$ is called a {\it standard sequence}. 
With the type of relative orientations of the sides, a 5-gon with outer connections is represented by a sequence uniquely. 
There are $4!=24$ standard sequences on each type, and we remark that there are some multiplicity by symmetries as a 5-gon of a spherical curve.

\phantom{x} 
\noindent {\bf Type 1:} 
A 5-gon of type 1 has two symmetries $\rho$ and $\varphi _1$. 
Two cyclic sequences which can be transformed into each other by some $\rho$s represent the same 5-gon with outer connections. 
For example, $abced$ and $aebcd$ represent the same 5-gon because $\rho (abced)=bcdae=aebcd$. 
Since the orientation is fixed, two sequences represent the same 5-gon when they are transformed into each other by a single $\varphi _1$ and orientation reversing (denoted by $\gamma$). 
For example, $abced$ and $acbde$ represent the same 5-gon because $\gamma ( \varphi _1 (abced))= \gamma (aedbc)=cbdea=acbde$. 
There are 8 equivalent classes of standard sequences up to some $\rho$s and a pair of $\varphi _1$ and $\gamma$: \\
\phantom{x} 

\noindent $abcde,$
$abced=abdce=acbde=acdeb=aebcd, $\\
$abdec=abecd=acdba=adbce=adebc, $\\
$abedc=adcbe=adecb=aecdb=aedbc, $\\
$acebd, $
$acedb=acbed=adceb=aebdc=aecbd, $
$adbec, $
$aedcb.$\\
\phantom{x} 
\begin{figure}[ht]
\begin{center}
\includegraphics[width=140mm]{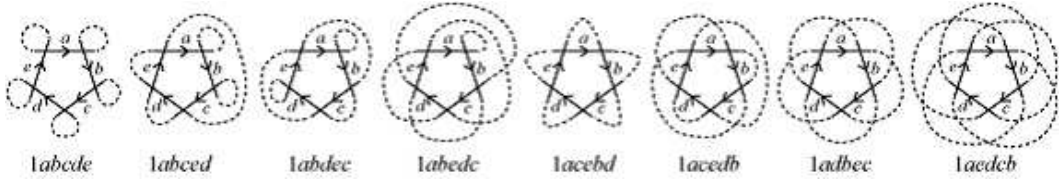}
\caption{The 5-gons of type 1. }
\label{type1}
\end{center}
\end{figure}

\noindent {\bf Type 2:} 
A 5-gon of type 2 has the reflection symmetry $\varphi _2$. 
Two cyclic sequences represent the same 5-gon when they are transformed into each other by a single $\varphi _2$ and orientation reversing $\gamma$. 
There are 16 equivalent classes of standard sequences up to a pair of $\varphi _2$ and $\gamma$: \\
\phantom{x}

\noindent $abcde, abced=aebcd, abdce=acdeb, abdec=acdbe, abecd, $\\
$abedc=aecdb, acbde, acbed=aecbd, acebd, acedb=aebdc, $\\
$adbce=adebc, adbec, adcbe=adecb, adceb, aedbc, aedcb $. \\
\phantom{x} 
\begin{figure}[ht]
\begin{center}
\includegraphics[width=140mm]{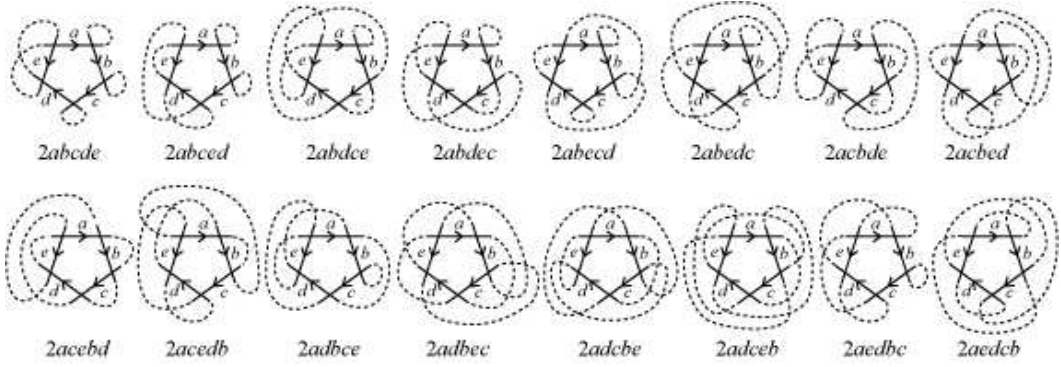}
\caption{The 5-gons of type 2. }
\label{type2}
\end{center}
\end{figure}

\noindent {\bf Type 3:} 
Two cyclic sequences represent the same 5-gon when they are transformed into each other by a single $\varphi _3$ and orientation reversing $\gamma$. 
There are 16 equivalent classes of standard sequences up to a pair of $\varphi _3$ and $\gamma$: \\
\phantom{x}

\noindent $ abcde, abced, abdce=aebcd, abdec=adebc, abecd=adbce,$\\
$ abedc=aedbc, acbde=acdeb, acbed=acedb, acdbe, acebd, adbec, $\\
$adcbe=aecdb, adceb=aecbd, adecb, aebdc, aedcb $. \\
\phantom{x} 
\begin{figure}[ht]
\begin{center}
\includegraphics[width=140mm]{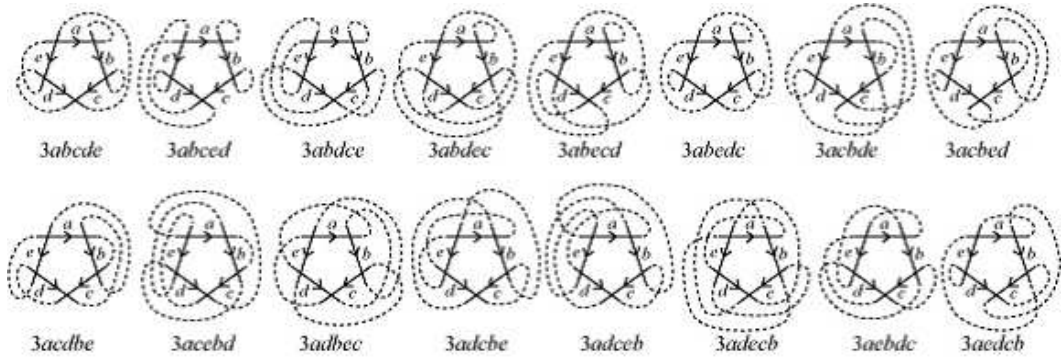}
\caption{The 5-gons of type 3. }
\label{type3}
\end{center}
\end{figure}

\noindent {\bf Type 4:} 
Two cyclic sequences represent the same 5-gon when they are transformed into each other by a single $\varphi _4$ and orientation reversing $\gamma$. 
There are 16 equivalent classes of standard sequences up to a pair of $\varphi _4$ and $\gamma$: \\
\phantom{x}

\noindent $abcde, abced=abdce, abdec=abecd, abedc, acbde=aebcd, $\\
$acbed=aebdc, acdbe=adebc, acdeb, acebd, acedb=adceb, $\\
$adbce, adbec, adcbe=aedbc, adecb=aecdb, aecbd, aedcb$. \\
\phantom{x} 
\begin{figure}[ht]
\begin{center}
\includegraphics[width=140mm]{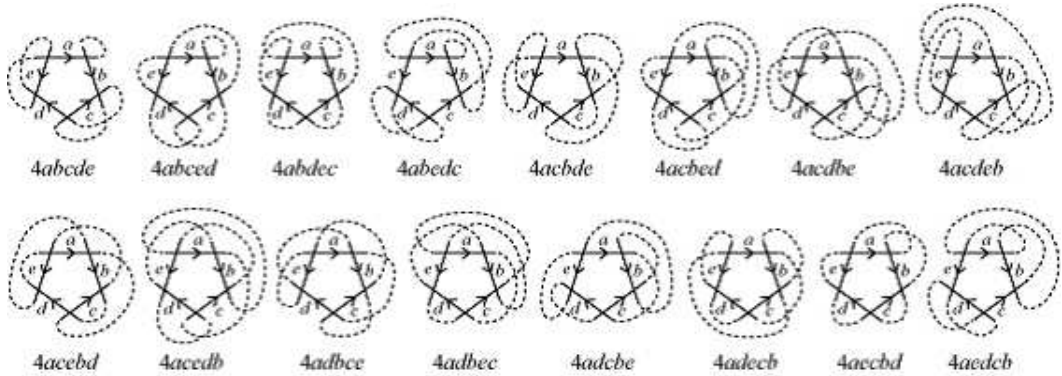}
\caption{The 5-gons of type 4. }
\label{type4}
\end{center}
\end{figure}

\noindent Thus, 5-gons are classified into the 56 types shown in Figure \ref{a-5-gons}. 
\hfill$\square$

\section{Proof of Theorem \ref{thm-u}}

In this section, Theorem \ref{thm-u} is proved. \\

\begin{figure}[ht]
\begin{center}
\includegraphics[width=25mm]{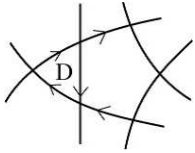}
\caption{Type 2 and 4. }
\label{type24}
\end{center}
\end{figure}

\begin{figure}[ht]
\begin{center}
\includegraphics[width=100mm]{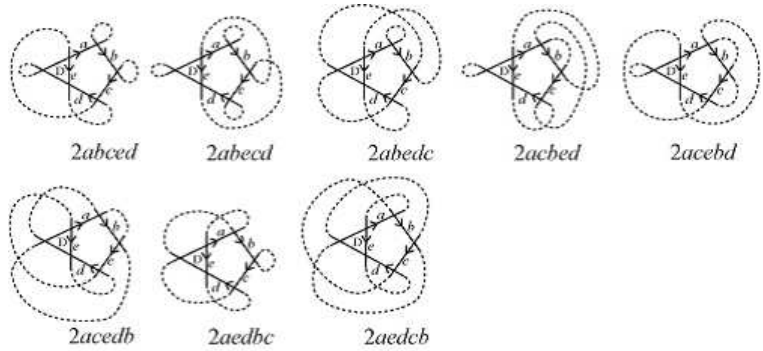}
\caption{The case that a 3-gon of type D and a 5-gon of type 2 share the side $e$. }
\label{u-type2}
\end{center}
\end{figure}

\begin{figure}[ht]
\begin{center}
\includegraphics[width=100mm]{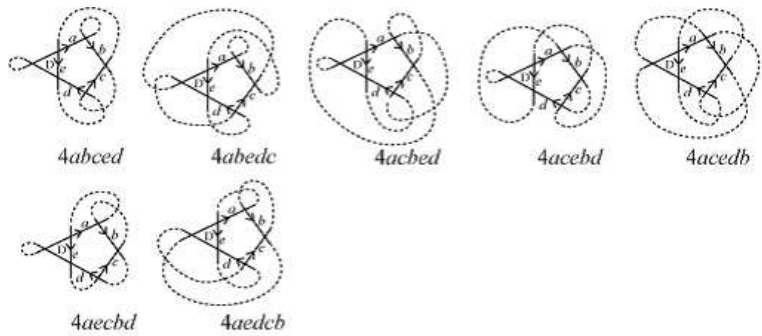}
\caption{The case that a 3-gon of type D and a 5-gon of type 4 share the side $e$. }
\label{u-type4-1}
\end{center}
\end{figure}

\begin{figure}[ht]
\begin{center}
\includegraphics[width=100mm]{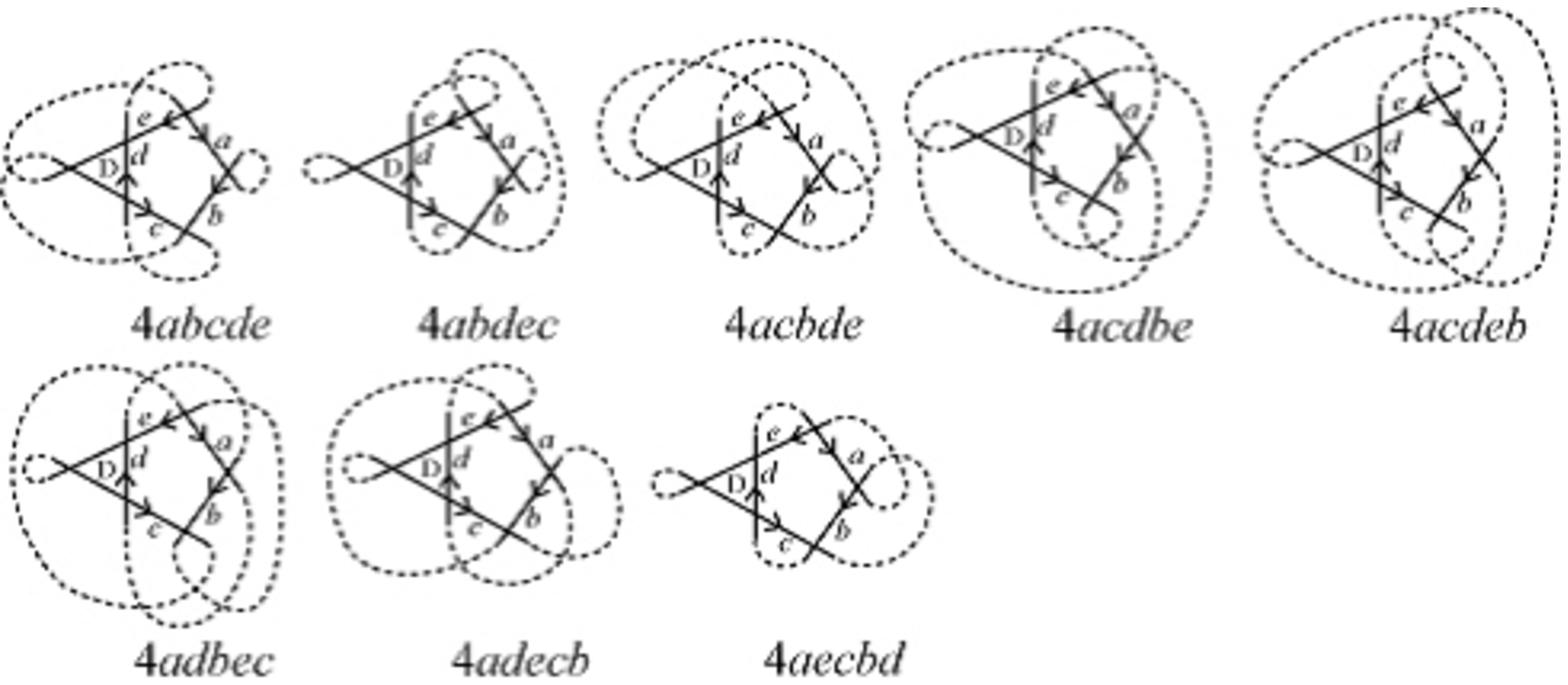}
\caption{The case that a 3-gon of type D and a 5-gon of type 4 share the side $d$. }
\label{u-type4-2}
\end{center}
\end{figure}

\begin{figure}[ht]
\begin{center}
\includegraphics[width=100mm]{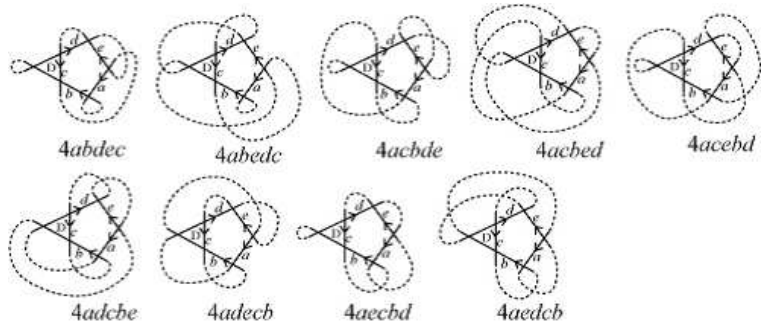}
\caption{The case a 3-gon of type D and a 5-gon of type 4 share the side $c$. }
\label{u-type4-3}
\end{center}
\end{figure}

\phantom{x}
\noindent {\rm Proof of Theorem \ref{thm-u}.} \ 
Let $P$ be a spherical curve with reductivity four. 
Since $P$ is reduced, the set $T_2$ in Figure \ref{u-set-r} is also an unavoidable set for $P$. 
Here, $P$ can not have the first and second configuration because they make reductivity three or less as discussed in \cite{S} and \cite{OS}. 
The third one of $T_2$ has already been discussed in Theorem 1 in \cite{OS}. 
Hence just the fourth one needs to be discussed here. 
Since the 3-gon should be of type D because 3-gons of type A, B and C make reductivity three or less, the 5-gon should be of type 2 or 4 with respect to the relative orientations of the sides (see Figure \ref{type24}). 
Let $a, b, c, d$ and $e$ be the sides of a 5-gon of type 2 and 4 as same as Figures \ref{type2} and \ref{type4}. 
When the 5-gon is of type 2, only the side $e$ can be shared with the 3-gon. 
In this case, by considering the outer connections of the 3-gon of type D, the 5-gon should be the one whose sequence includes $a$, $e$, $d$ with this cyclic order, which are the type of $2abced$, $2abecd$, $2abedc$, $2acbed$, $2acebd$, $2acedb$, $2aedbc$ and $2aedcb$. 
Hence the eight configurations with outer connections illustrated in Figure \ref{u-type2} are obtained. 
When the 5-gon is of type 4, the sides $e$, $d$ and $c$ can be shared with the 3-gon. 
When $e$ is shared, the 5-gon should be the one whose sequence includes $a$, $e$, $d$ with this cyclic order, which are the type of $4abced$, $4abedc$, $4acbed$, $4acebd$, $4acedb$, $4aecbd$ and $4aedcb$. 
Hence the seven configurations with outer connections in Figure \ref{u-type4-1} are obtained. 
When $d$ is shared, the 5-gon should be the one whose sequence includes $c$, $d$, $e$ with this cyclic order, which are the type of $4abcde$, $4abdec$, $4acbde$, $4acdbe$, $4acdeb$, $4adbec$, $4adecb$ and $4aecbd$. 
Hence the eight configurations with outer connections in Figure \ref{u-type4-2} are obtained. 
When $c$ is shared, the 5-gon should be the one whose sequence includes $b$, $d$, $c$ with this cyclic order, which are the type of $4abdec$, $4abedc$, $4acbde$, $4acbed$, $4acebd$, $4adcbe$, $4adecb$, $4aecbd$ and $4aedcb$. 
Hence the nine configurations with outer connections in Figure \ref{u-type4-3} are obtained. 

\noindent Hence, the set $U_3$ is an unavoidable set for a spherical curve of reductivity four. \\

\hfill$\square$

\section{Appendix: 5-gons on chord diagrams}

\begin{figure}[ht]
\begin{center}
\includegraphics[width=140mm]{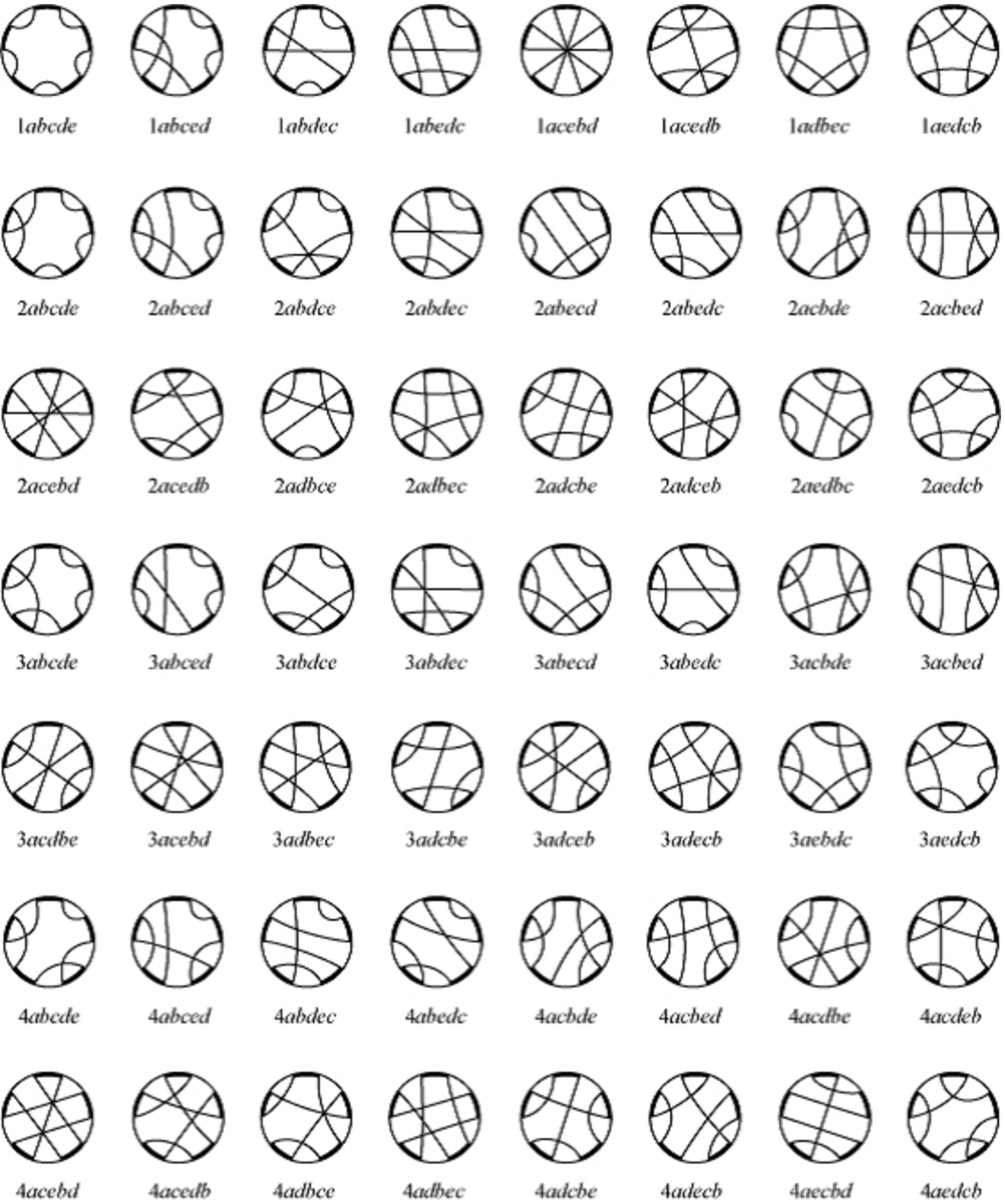}
\caption{All the 5-gons on chord diagrams. There are no endpoints of segments in the interior of each thick arc. }
\label{a-5-gons-co}
\end{center}
\end{figure}

A {\it chord diagram} of a spherical curve $P$ is a preimage of $P$ with each pair of points corresponding to the same double point connected by a segment as $P$ is assumed to be an image of an immersion of a circle to $S^2$. 
In Figure \ref{a-5-gons-co}, all the 5-gons of a spherical curve on chord diagrams are listed. \\


\phantom{x}
\noindent {\rm ACKNOWLEDGMENTS.} \ 
The authors are grateful to the members of COmbinatoric MAthematics SEMInar (COMA Semi) for helpful comments. 
They also thank the timely help given by Yuki Miyajima in discovering reduced spherical curves without 2-gons and 3-gons of type A and B. 
The second author was supported by Grant for Basic Science Research Projects from The Sumitomo Foundation (160154).

\phantom{x}

\noindent Department of General Systems Studies, \\
\noindent University of Tokyo, \\
\noindent 3-8-1, Komaba, Meguro-ku, Tokyo 153-8902, Japan. \\
\noindent Email: kashiwa@idea.c.u-tokyo.ac.jp\\

\noindent Department of Mathematics, \\
\noindent National Institute of Technology, Gunma College, \\
\noindent 580 Toriba-cho, Maebashi-shi, Gunma 371-8530, Japan. \\
\noindent Email: shimizu@nat.gunma-ct.ac.jp

\end{document}